\documentclass[review]{elsarticle}

\usepackage{lineno,hyperref}
\modulolinenumbers[1]

\journal{arXiv}

 \biboptions{numbers,sort&compress}
  
\usepackage{graphicx}
\usepackage{epstopdf}
\usepackage{amsmath}
\usepackage{amsfonts}
\usepackage[compact]{titlesec}
\usepackage{amsthm}
\usepackage{amssymb}
\usepackage{multirow}
\usepackage{color}
\usepackage{mathtools}
\usepackage{geometry}
\usepackage{times}
\usepackage{algorithm}
\usepackage{algorithmicx, algpseudocode}

\usepackage{verbatim}

\theoremstyle{definition}

\definecolor{Green}{rgb}{0,.5,0}
\definecolor{Blue}{rgb}{0,.1,.85}
\definecolor{Cyan}{rgb}{.2,.6,.7}
\definecolor{Purple}{rgb}{.5,0,1}
\definecolor{deepred}{rgb}{.8,.1,.2}
\newcommand{\rev}[1]{\textcolor{black}{#1}}

\begin{document}

\begin{frontmatter}

\title{Soft \rev{IsoGeometric A}nalysis of the Bound States of a Quantum Three-Body Problem in 1D} 

\author[ad]{Danyang Li\corref{corr}}
\ead{Danyang.Li@anu.edu.au}

\author[ad]{Quanling Deng}
\cortext[corr]{Corresponding author}
\ead{Quanling.Deng@anu.edu.au}

\address[ad]{School of Computing, Australian National University, Canberra, ACT 2601, Australia}

\begin{abstract}
The study of quantum three-body problems has been centred on low-energy states 
that rely on accurate numerical approximation. 
Recently, isogeometric analysis (IGA) has been adopted to solve the problem 
as an alternative but more robust (with respect to atom mass ratios) method 
that outperforms the classical Born-Oppenheimer (BO) approximation\rev{, especially for the cases with small mass ratios}.
In this paper, we focus on the performance of IGA 
and apply the recently-developed softIGA to \rev{further reduce} the spectral errors of the low-energy bound states. 
This is an extension to the recent work that is published as an ICCS conference paper in  \cite{deng2022isogeometric}.
The main idea of softIGA is to add high-order derivative-jump terms with a penalty parameter to the IGA bilinear forms. 
With an optimal choice of the penalty parameter, 
we observe eigenvalue error superconvergence.
\rev{Herein, the optimal parameter coincides with the ones for the Laplace operator (zero potential) and can be heuristically computed for a general elliptic operator.}
We focus on linear (finite elements) and quadratic elements and demonstrate the outperformance of softIGA over IGA through a variety of examples including both two- and three-body problems in 1D.
\end{abstract}


\end{frontmatter}

\section{Introduction}

The dynamics of three interacting bodies constitute one of the classical challenges in physics and contain many unsolved questions. 
\rev{The three-body problem arises in application in quantum mechanics that models the motion of three particles \cite{schmid2017quantum,eyges1959quantum}.} 
The main challenges remain in solving the underlying Schr\"odinger equation. 
To reduce the computational cost, 
the first idea is to reduce the dimensionality by using the center-of-mass of the many-body system.
This is also referred to as the hyperspherical approach \cite{raynal1970transformation}.
For example, for a system of $n$-body problems in 1D, 
the overall reduced dimension is $n-1$ after removing the center-of-mass motion. 
In this paper, we focus on two- and three-body problems in 1D.
Thus, the reduced dimension is one and two, respectively. 

The three-body problem generally does not have \rev{analytical solution}. Thus, one solves the Schr\"odinger equation numerically \cite{berezin1991schrodinger}.  
One of the commonly used approximation methods when considering the solution of quantum mechanical equations for systems involving electrons and nuclei is the Born-Oppenheimer (BO) approximation (also known as an adiabatic approximation) \cite{born1985quantentheorie}. 
In molecular dynamics, one may consider the nuclei and their surrounding electrons as a many-body system. 
In general, the electrons are much lighter than nuclei. 
This allows the wave functions or the solution states of atomic nuclei and electrons in a molecule to be treated separately.
For example, \cite{cederbaum2008born} discussed the nuclear dynamics in the framework of 
a fully time-dependent BO approximation. 
BO was also adopted in \cite{happ2019universality} to establish the universality in 
a one-space dimensional heavy-heavy-light three-body system. 
Except for BO approximation, other methods such as explicitly correlated Gaussians (ECGs)  
can be applied to solve few-body systems \cite{mitroy2013theory}. 
The Skorniakov and Ter-Martirosian (STM) 
method was developed to solve the three-body bound states in the limiting case of zero-range forces \cite{skorniakov1957three}.

However, current numerical methods are unsatisfactory in terms of accuracy, robustness, and reliability.
These aspects are partially addressed in the recent work \cite{deng2022isogeometric}
by applying the \rev{isogeometric analysis (IGA)} tools developed in \cite{cottrell2009isogeometric,hughes2005isogeometric}.
In this paper, \rev{we extend the recent work \cite{deng2022isogeometric} that was published as an ICCS conference
paper and}
 develop an advanced IGA method to solve the low-energy bound states of the one-dimensional quantum two- and three-body problems. 
 \rev{The major new contribution is on the application of the new method and the demonstration of its outperformance by comparing it with the one in \cite{deng2022isogeometric}.
 }
The method is robust in the sense of arbitrary mass ratios 
and any interaction potentials as long as they lead to bound states.
Based on our previous studies \cite{deng2022isogeometric} of finite element analysis (FEA) methods 
and more advanced methods of \rev{IGA}, 
we introduce the soft isogeometric analysis (SoftIGA) method for this purpose. 
SoftIGA, presented in \cite{deng2023softiga}, mainly extend the idea of softFEM developed in \cite{deng2021softfem}
to the IGA setting. 
For $p$-th order IGA elements with maximal continuity, the basis functions are $C^{p-1}$-continuous. 
The jumps appear when taking the basis functions’ $p$-th order partial derivatives. 
We thus penalize this $p$-th order derivative-jump 
and subtract from IGA (for $p=1,2$; outlier-free IGA \cite{hiemstra2021removal,manni2022application,deng2021outlier} for $p\ge 3$) bilinear form  an inner product of the derivative-jumps 
of the basis functions in both trial and test spaces. 
The previous result shows that softIGA is able to reduce the stiffness 
(consequently, the condition numbers) of the IGA discretized problem. 
Therefore, we adopt softIGA to solve the quantum three-body problem as a second-order differential eigenvalue problem. 
\rev{Finite element method (FEM) and softFEM are also applied to solve the problem and results are compared.}

We organise the rest of the paper as follows. In Section \ref{sec:ps}, we state the differential eigenvalue problem that unifies the two- and three-body problems in 1D. 
In section \ref{sec:m}, we first review the IGA discretization method and then describe the softIGA discretization method in subsection \ref{sec:softIGA}.
In particular, when using linear elements, softIGA reduces to softFEM.
Section \ref{sec:num} collects and discusses various numerical tests to demonstrate the performance of the proposed method. 
Concluding remarks are presented in section \ref{sec:conc}.

\section{Problem Statement} \label{sec:ps}

The heavy-light two-body problems and heavy-heavy-light three-body problems are modelled as the dimensionless stationary Schr\"odinger equations \cite{happ2019universality, happ2022universality}. As in  \cite{deng2022isogeometric}, we generalize these problems to any mass ratio and unify it into the differential eigenvalue problem of finding the eigenpair $(\lambda, u)$ such that
\begin{equation} \label{eq:schro}
   -\nabla \cdot (\kappa\nabla u) - \gamma u = \lambda u\qquad\forall x\in \Omega
\end{equation}
where $\nabla$ is the gradient operator, $\nabla\cdot$ is the divergence operator, and $u$ denotes an eigenstate. 
\rev{ $\Omega$ is an infinite domain, $\Omega = \mathbb{R}$ for two-body problem and $\Omega = \mathbb{R}^2$ for three-body problem.}
Herein, the differential operator is also referred to as the Hamiltonian.
The two-body problem is the unified problem \eqref{eq:schro} in 1D while the three-body problem is \eqref{eq:schro} in 2D.
The potential function is defined as
\begin{equation}
\gamma = 
\begin{cases}
v(x), \quad  & \text{in 1D}, \\
v(x+y/2) + v(x-y/2), \quad  & \text{in 2D}, 
\end{cases}
\end{equation}
where 
$
    v(\xi) = \beta f(\xi),
$
  $\beta$ denotes a magnitude, and $f$ denotes the shape of the interaction potential. 
The diffusion coefficient is defined as 
\begin{equation}
\kappa = 
\begin{cases}
\frac12, \quad  & \text{in 1D}, \\
(\frac{\alpha_x}{2},0;0,\frac{\alpha_y}{2}) \ \text{being a tensor matrix}, \quad  & \text{in 2D}.
\end{cases}
\end{equation}
\rev{$\alpha_x$ and $\alpha_y$ are the coefficient which related to heavy particles $m_h$ and light particles $m_l$. They are defined as}
\begin{equation}
\alpha_x = \frac{1/2 + m_h/m_l}{1 + m_h/m_l} \quad \alpha_y = \frac{2}{1 + m_h/m_l}
\end{equation}

The direct finite element discretization of the differential operator $\mathcal{L} = -\nabla\cdot(\kappa\nabla)-\gamma$ 
leads to a stiffness matrix 
which is not necessarily positive and definite.
This in turn introduces a potential problem when solving the resulting linear algebra problem.
This means that some eigenvalues are negative,  
corresponding to the case that the attractive interactions lead to negative eigenvalues. 
To avoid negative eigenvalues at the discretization level, 
we rewrite $-\nabla \cdot (\kappa\nabla u)-\gamma u = \lambda u$ by adding a positive constant $\gamma_0$
 to obtain $-\nabla \cdot (\kappa\nabla u)-(\gamma-\gamma_0) u = (\lambda-\gamma_0) u$ 
 such that $\gamma-\gamma_0 < 0$ for all $x \in \Omega$. With this in mind, problem \eqref{eq:schro} can be rewritten as
\begin{equation} \label{eq:pde}
  -\nabla \cdot (\kappa\nabla u)+\rev{\hat{\gamma}} u = \lambda u\qquad\forall x\in \Omega,
\end{equation}
where \rev{$\hat{\gamma} = \gamma_0 -\gamma$. }

In quantum mechanics, when $f$ is symmetric and describes a short-range interaction,
the lower-energy eigenfunctions are localized and bounded in a small subdomain
where the higher-energy eigenfunctions are scattered to the overall infinite domain.
Mathematically, this requires that $|\xi|^2f(|\xi|) \rightarrow 0$ as $|\xi| \rightarrow \infty$.

\section{Soft Isogeometric Analysis} \label{sec:m}

In this section, we first review the standard IGA discretization, then present softIGA method for the unified problem \eqref{eq:pde} on a finite domain $\Omega_\epsilon = [-x_\epsilon,x_\epsilon]^d$, $d = 1$, $2,$ with homogeneous boundary condition

\begin{equation} \label{eq:bound_cond}
   u = 0,\qquad\forall x\in \partial\Omega_\epsilon.
\end{equation}
Herein, we focus on finding the bound states of \eqref{eq:pde} that their function values vanish at $\pm \infty.$
Hence, $\Omega_\epsilon$ is viewed as an approximation of the domain $\Omega = \mathbb{R}^d$ for \eqref{eq:pde}.

\subsection{Variational Formulation}

\rev{Assuming $\Omega_\epsilon=[-x_\epsilon,x_\epsilon]^d \subset\mathbb{R}^d$, $d =$ 1,2 is a bounded domain Lipschitz boundary $\partial\Omega_\epsilon$. We adopt the standard notation for the Hilbert and Sobolev spaces. 
We denote by $(\cdot, \cdot)_{\Omega_\epsilon}$ and $\|\cdot\|_{\Omega_\epsilon}$ the $L^2$-inner product and its norm, respectively. 
}
\rev{Let $H_0^1(\Omega_{\epsilon})$ be the Sobolev space with functions in $H^1(\Omega_{\epsilon})$ that are vanishing at the boundary.}

The variational formulation of \eqref{eq:pde} at the continuous level is to find the eigenvalue $\lambda \in\rev{\mathbb{R}}$ and the associated eigenfunction \rev{$u \in H_0^1(\Omega_{\epsilon})$ with $\|u\|_{\Omega_{\epsilon}}=1$ such that}

\begin{equation}\label{eq:var_formul}
   a(w,u) = \lambda b(w,u),\qquad\forall w\in H_0^1(\Omega_\epsilon),
\end{equation}
where the bilinear forms are defined as for $v,w \in H_0^1(\Omega_\epsilon)$

\begin{equation}\label{eq:bil_form}
   a(w,u) := (\kappa\nabla v, \nabla w)_{\Omega_\epsilon} + (\rev{\hat{\gamma}} v,w)_{\Omega_\epsilon},\qquad b(w,u) := (v,w)_{\Omega_\epsilon}.      
\end{equation}

The eigenvalue problem \eqref{eq:var_formul} with $\rev{\hat{\gamma}} > 0$ in domain $\Omega_\epsilon$, has
a countable set of positive eigenvalues $\{\lambda_j\}_{j=1}^\infty$
and an associated set of orthonormal eigenfunctions $\{u_j\}_{j=1}^\infty$, meaning, $(u_j,u_k) = \delta_{jk}$, where $\delta_{jk} = 1$ is the Kronecker delta. Furthermore, the eigenfunctions are also orthogonal in the energy inner product since $a(u_j,u_k)=\lambda_jb(u_j,u_k)=\lambda_j\delta_{jk}$.
For bound states, this set of eigenpairs is finite.

\subsection{Isogeometric Analysis} \label{sec:iga}

Standard IGA adopts the Galerkin finite element analysis framework at the discrete level.
We first divide the boundary domain $\Omega_\epsilon$ with a uniform tensor-product mesh. 
\rev{Let $\tau$ and $\mathcal{T}_h$ be a general element and its collection such that $\bar{\Omega}_\epsilon=\cup_{\tau \in \mathcal{T}_h}\tau$. Let $h = \max_{\tau \in \mathcal{T}_h}$ diameter($\tau$).}
The isogeometric analysis of \eqref{eq:pde} in variational formulation seeks $\lambda^h \in \mathbb{R}$ and $u^h \in V_p^h$ with $\|u^h\|_{\Omega_\epsilon} = 1$ such that

\begin{equation}
    a(w^h,u^h)= \lambda^h b(w^h,u^h),\qquad\forall w^h\in V_p^h,
\end{equation}
where $V_p^h \subset H_0^1(\Omega_\epsilon)$ donates the IGA approximation space
that is to be specified below.

The IGA approximation space consists of B-splines as basis functions. 
We construct using the Cox-de Boor recursive formula in 1D. Let $X = \{x_0,x_1,\cdots,x_m\}$ be a knot vector with \rev{non-decreasing sequential knots $x_j$}. The \rev{$j$-th B-spline basis function of degree $p$ for the space $V^h_p$}, denoted as $\phi_p^j(x)$, is defined recursively as (see also in \cite{de1978practical})
\begin{equation} \label{eq:Bspline}
\begin{aligned}
\phi^j_0(x) & = 
\begin{cases}
1, \quad \text{if} \ x_j \le x < x_{j+1}, \\
0, \quad \text{otherwise}, \\
\end{cases} \\ 
\phi^j_p(x) & = \frac{x - x_j}{x_{j+p} - x_j} \phi^j_{p-1}(x) + \frac{x_{j+p+1} - x}{x_{j+p+1} - x_{j+1}} \phi^{j+1}_{p-1}(x).
\end{aligned}
\end{equation}

With this in mind, the IGA approximation space in 1D is $V_p^h= \text{span} \{\phi_p^j\}_{j=1}^{N_h}$. 
In 2D, based on the tensor-product structure, $V_p^h= \text{span} \{\phi_{p_x}^{j_x}(x)\phi_{p_y}^{j_y}(y)\}_{j_x,j_y=1}^{N_x,N_y}$
where $p_x$, $p_y$ specify the approximation order in each dimension. $N_x$ and $N_y$ are the total numbers of basis functions in each dimension and $N_h$ is the total number of degrees of freedom. 
\rev{Herein, for repeated knots at the boundary nodes, we apply the end-node limits; we refer to \cite{buffa2011isogeometric,evans2013isogeometric} for more details.} 
Throughout the rest of the paper, we focus on $C^0$ linear finite elements and $C^1$ quadratic IGA elements \rev{(using maximal continuity of the B-splines)}.




\subsection{SoftIGA Discretization} \label{sec:softIGA}

Following the IGA variational formulation, 
the softIGA is to find $\tilde{\lambda}^h \in \rev{\mathbb{R}}$ and $\tilde{u}^h \in V_p^h$ with \rev{$\|\tilde{u}^h\|_{\Omega_{\epsilon}} = 1$ }such that

\begin{equation} \label{eq:vfhh}
  a(w^h,\tilde{u}^h) -\eta s(w^h,\tilde{u}^h)= \tilde{\lambda}^h b(w^h,\tilde{u}^h),\qquad\forall w^h\in V_p^h,
\end{equation}
where $\eta \in [0, \eta_{\max})$ is the softness parameter, $s(w,v)$ is the softness bilinear form defined over the interfacial derivative jumps to be specified below.

For a tensor-product mesh $\mathcal{T}_h$, 
let $F$ denote a face while $\mathcal{F}$ represents the set of interior faces of the mesh,
and let $F_b$ denote a face at the boundary while $\mathcal{F}_b$ denote the set of boundary faces of the mesh. 
As in the softFEM, we introduce the $p$-th order derivative-jump as:

\begin{equation}
\begin{aligned}
  [\tilde{\nabla}^p v\cdot n] = \tilde{\nabla}^p v\vert_{F_1}\cdot n_1 + \tilde{\nabla}^p v\vert_{F_2}\cdot n_2, \qquad & \forall v \in V_p^h, \\
  [\tilde{\nabla}^pv\cdot n_{F_b}] = \tilde{\nabla}^pv|_{E}\cdot n_{F_b}, \qquad & F_b \in \mathcal{F},
\end{aligned}
\end{equation}
where $\tilde{\nabla}^p=(\partial_{x_1}^p,\cdots,\partial_{x_d}^p)^T$. 
Herein, for any two neighbouring elements $\tau_1 $ and $\tau_2$ with an interface $F =  \tau_1 \cap \tau_2$, 
$F_1 := \tau_1 \cap F \in \mathcal{F}$ and $F_2 := \tau_2 \cap F \in \mathcal{F}$. $n_1, n_2,$ and $n_{F_b}$ are the outward unit normal vectors.
The softness bilinear form $s(w,v)$ is then defined as

\begin{equation}
    s(w,v) = \left\{ \begin{array}{ll} \sum_{F \in \mathcal{F}} h^{2p-1}([\kappa \tilde{\nabla}^pw\cdot n],[\tilde{\nabla}^pv\cdot n]),\qquad & p \text{ is odd,} \\
    \sum_{F \in \mathcal{F}} h^{2p-1}([\kappa \tilde{\nabla}^pw\cdot n],[\tilde{\nabla}^pv\cdot n])\\
    +2\sum_{F_b \in \mathcal{F}_b} h^{2p-1}(\kappa [\tilde{\nabla}^pw\cdot n],[\tilde{\nabla}^pv\cdot n]),\qquad & p \text{ is even.}
    \end{array}\right.
\end{equation}

We set $\eta \in [0,\eta_{max})$ where $\eta_{max}$ is to be determined such that $\tilde{a}(w^h,\tilde{u}^h) = a(w^h,\tilde{u}^h) - \eta s(w^h,\tilde{u}^h)$is coercive. We will discuss the value of $\eta_{max}$ through numerical experiments in the next section.
This softIGA formulation \eqref{eq:vfhh} leads to a generalized matrix eigenvalue problem

\begin{equation}\label{eq:gmevp}
    \tilde{K}\tilde{U} = \tilde{\lambda}^h M\tilde{U},\qquad \tilde{K}:= K-\eta S,
\end{equation}
where $K_{kl} = a(\phi_p^k,\phi_p^l)$, $M_{kl} = b(\phi_p^k, \phi_p^l)$, $S_{kl}=s(\phi_p^k,\phi_p^l)$ and $\tilde{U}$ is the eigenvector representing the coefficients of the B-spline basis functions.

\section{Numerical Experiments} \label{sec:num}

In this section, we provide various numerical examples to demonstrate the performance of softIGA for linear (also referred to as softFEM) and quadratic elements. 
In particular, we consider both polynomial and exponential decaying potentials that lead to \rev{bound state}. 
We also show numerical examples with mass ratios that are of different scales, namely $10^k, k=0,1,2$.  


\subsection{SoftIGA Discretization Accuracy}
We focus on the two- and three-body problems with a potential of polynomial decay (cube of a Lorentzian)
\begin{equation}\label{eq:poly_potential}
    f(\xi)=\frac{1}{(1+\xi^2)^3}
\end{equation}
and a potential of exponential decay (Gaussian)
\begin{equation}\label{eq:exp_potential}
    f(\xi)=e^{-\xi^2}.
\end{equation}

 For these potentials, it is impossible to find exact analytical solutions to the problem \eqref{eq:pde}. 
 To characterize the error, we use $C^6$ septic IGA elements with a fine mesh to solve the problem \eqref{eq:pde} for reference solutions. 
 We focus on the eigenvalue error that is defined as
\begin{equation}\label{eq:eig_err}
    e_j = |\tilde{\lambda}_j^h- \lambda_j| \approx |\tilde{\lambda}_j^h-\hat{\lambda}_j|,
\end{equation}
where $\tilde{\lambda}_j^h$ is a softIGA eigenvalue and $\hat{\lambda}_j$ is a reference eigenvalue that is of high (much higher than $\tilde{\lambda}_j^h$) accuracy approximating the exact one $\lambda_j$.

\begin{figure}
    \centering 
    \includegraphics[width=11cm]{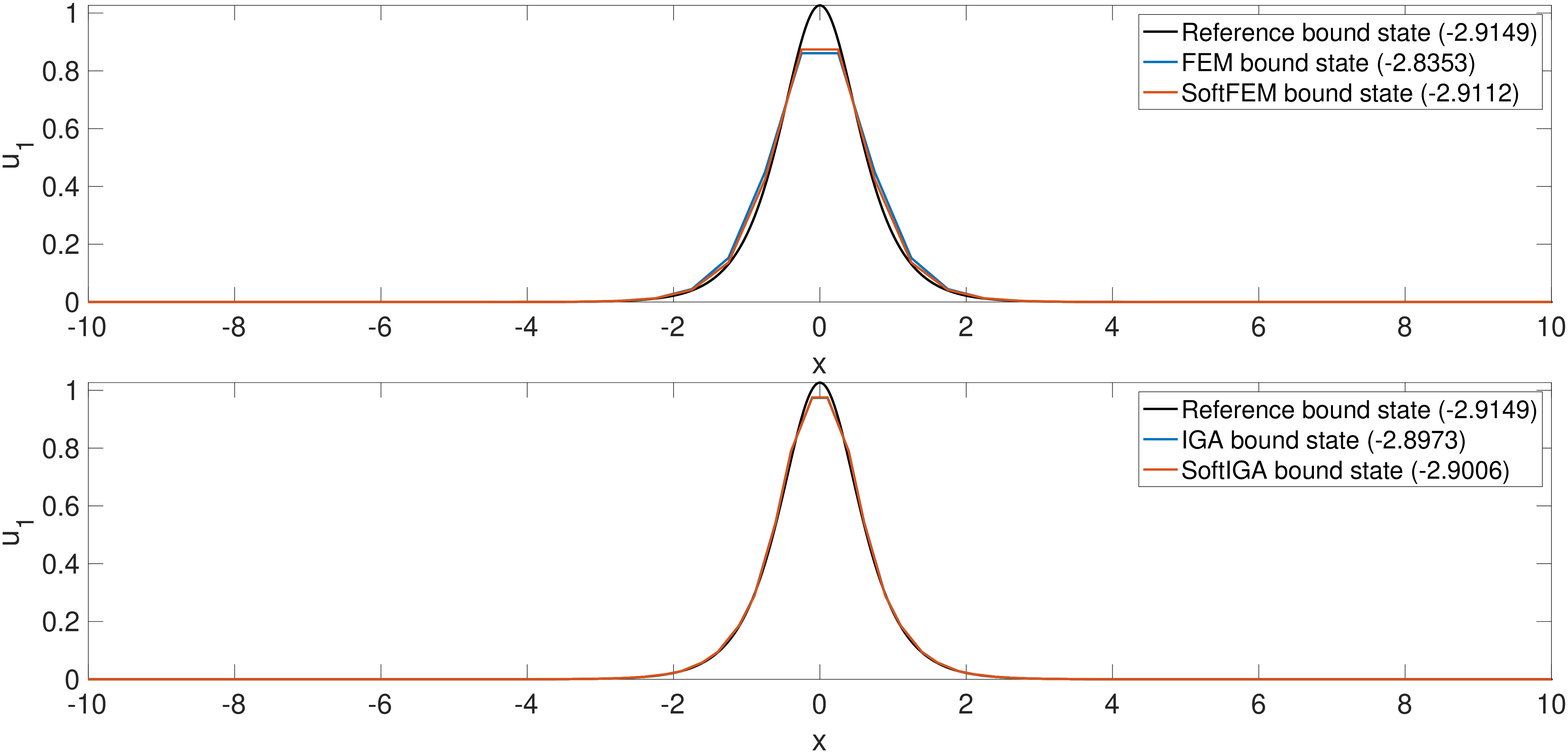}
    \caption{The first bound state of the two-body problem with a mass ratio 20, potential magnitude $\beta  = 5$, and \rev{a polynomial} decaying potential \eqref{eq:poly_potential} using $C^0$ linear softFEM and $C^1$ quadratic softIGA.}
    \label{fig:u1}
%
    \centering
    \includegraphics[width=11cm]{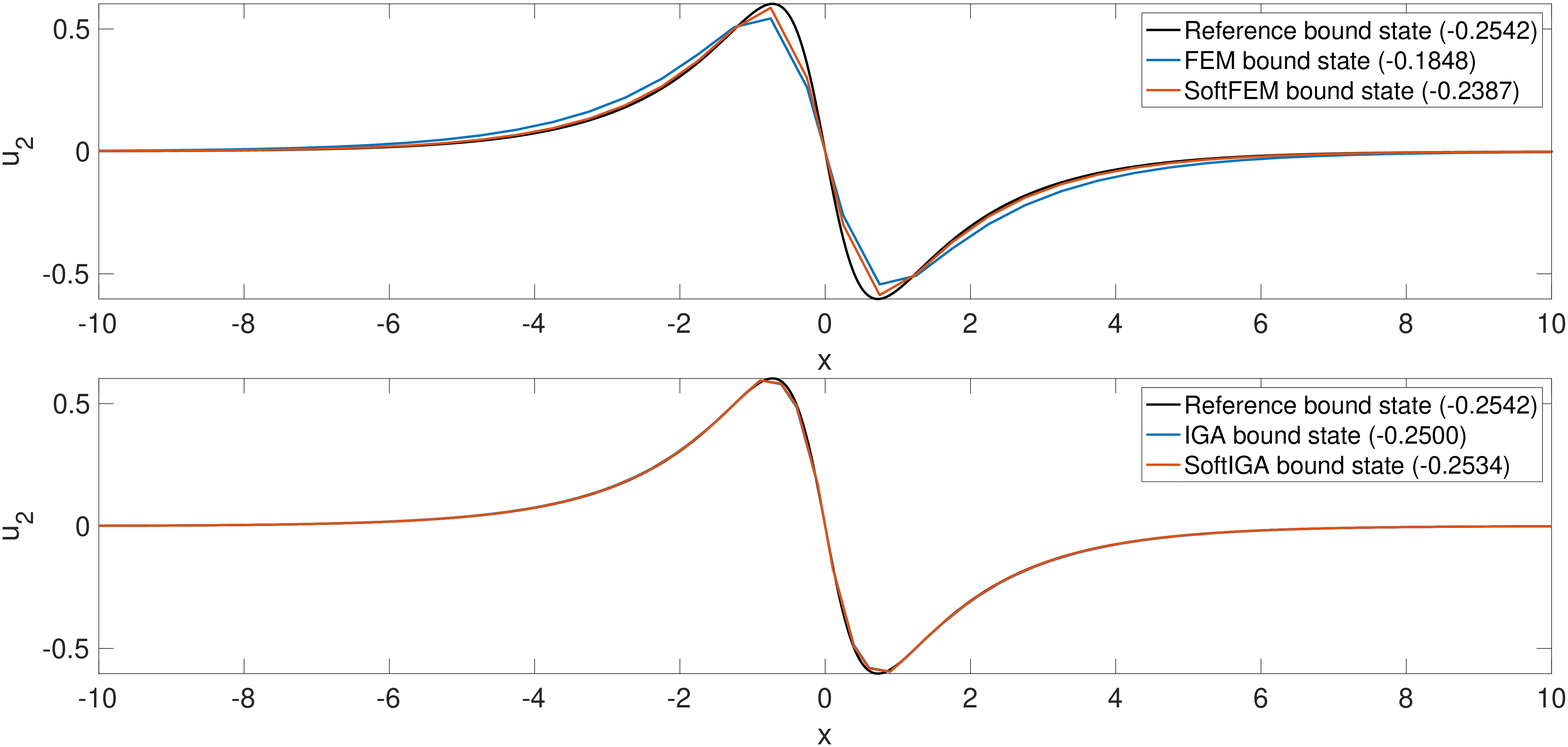}
    \caption{The second bound state of the two-body problem with a mass ratio 20, potential magnitude $\beta  = 5$, and \rev{a polynomial} decaying potential \eqref{eq:poly_potential} using $C^0$ linear softFEM and $C^1$ quadratic softIGA.}
    \label{fig:u2}
\end{figure}

We first consider the bound states with the polynomial decaying potential \eqref{eq:poly_potential} with magnitude $\beta = 5$ on domain $\Omega_\epsilon = [-20, 20]$.
Figures \ref{fig:u1} and \ref{fig:u2} show the two bound states approximated using both $C^0$ linear \rev{FEM/softFEM} and $C^1$ quadratic \rev{IGA/softIGA} with 80 uniform elements. 
\rev{For high-order accuracy, the reference eigenpair solutions are obtained by
using $C^6$ septic IGA (using $C^6$ septic  B-splines as basis functions) with 5000 uniform elements, }
which lead to the reference eigenvalue \rev{$\hat{\lambda}_1 = -2.9149185630$ (in the numerical implementation, a shift operator with $\gamma_0=5$ was used so that the discretised matrices are positive definite; see Section \ref{sec:ps})}. 
We observe that \rev{all the methods capture the
main behaviour of the bound states with a small number of elements (herein the mesh size is $0.5$). The ``soft" does not improve much on the accuracy of the eigenfunctions but significantly on the eigenvalues. 
We thus focus on the numerical study of the eigenvalues in this section.}

\begin{figure}
    \centering
    \includegraphics[width=\textwidth]{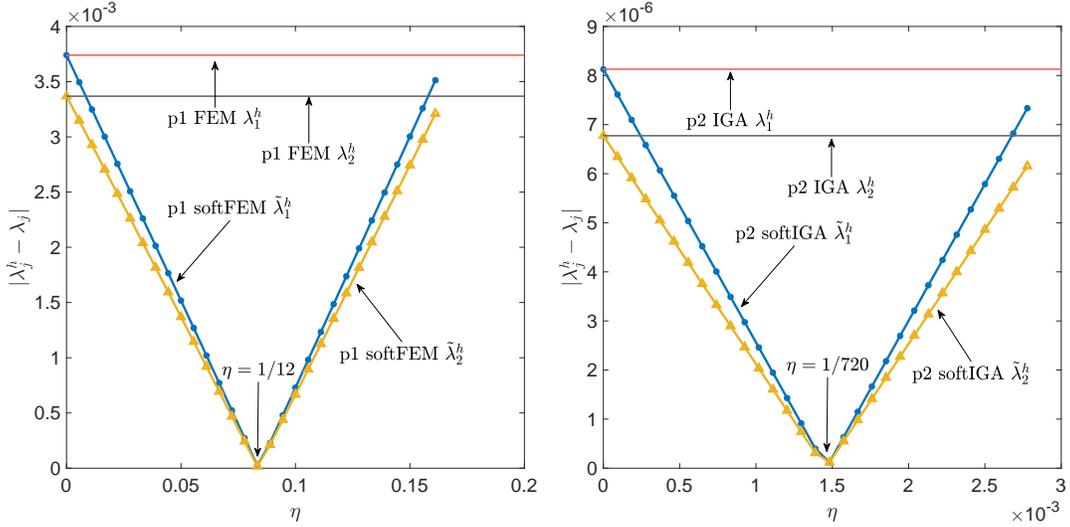}
    \caption{Eigenvalue errors of $C^0$ linear softFEM (left plot) and $C^1$ quadratic softIGA (right plot) and for the two-body problem with potential \eqref{eq:poly_potential} and $\beta = 5$ on domain $\Omega_\epsilon = [-20,20]$.}
    \label{fig:c0c1lambda1_2}
\end{figure}

\begin{figure}
    \centering
    \includegraphics[width=\textwidth]{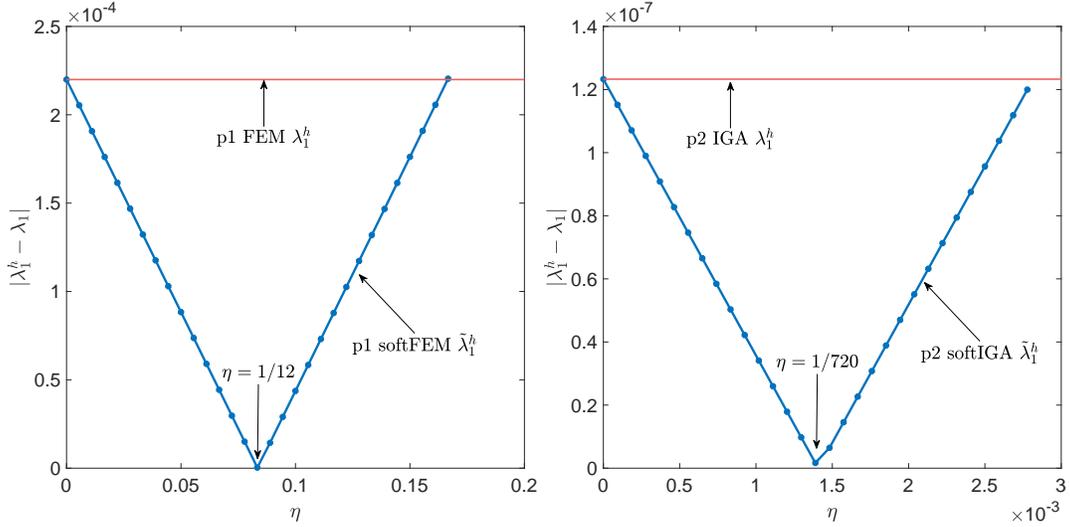}
    \caption{Eigenvalue errors of $C^0$ linear softFEM (left plot) and $C^1$ quadratic softIGA (right plot) and for the one-body problem with potential \eqref{eq:exp_potential} and $\beta = 1$ on domain $\Omega_\epsilon = [-20,20]$.}
    \label{fig:c0c1lambda1}
\end{figure}

The performance of softIGA depends on the softness parameter $\eta$. 
We now study the impact of $\eta$ on the approximation accuracy.
Figure \ref{fig:c0c1lambda1_2} shows the eigenvalue errors of  two bound states with respect to $\eta$.
The left plot shows the case of softFEM while the right plot shows that of the softIGA. 
The errors for the second eigenvalue are smaller than the first one.
\rev{This is due to the fact that} for these bound states, the eigenvalues are negative and $|\lambda_2| < |\lambda_1|.$
The parameter and domain size setting is the same as in Figure \ref{fig:u1}. 
We use 400 uniform elements in both cases.
To characterize the eigenvalue error, we use as an accurate approximation a reference eigenvalue solved by $C^6$ septic IGA with 5000 uniform elements.
We observe that for $\eta \in (0, 1/6)$, the eigenvalue errors of softFEM are smaller than those of FEM.
Similarly, for $\eta \in (0, 1/360)$, the eigenvalue errors of softIGA are smaller than those of IGA. 
For $C^0$ linear softFEM, the eigenvalue error is minimized when $\eta=\frac{1}{12}$. The errors are $2.21\times 10^{-5}$ and $1.61\times 10^{-5}$
for the first and second eigenvalue, respectively. 
\rev{We also observe that the eigenvalue errors have approximately a linear dependence on the distance of $\eta$ to the optimal one.}
Similarly, the eigenvalue error of softIGA is minimized when $\eta = \frac{1}{720}$.
The errors are $1.17\times 10^{-7}$ and $1.20\times 10^{-7}$ for the first and second eigenvalue, respectively.
These optimal values (and their ranges) for the softness parameter $\eta$ are matching with the ones in \cite{deng2021softfem} for softFEM and in \cite{deng2023softiga} for softIGA. 
Lastly, Figure \ref{fig:c0c1lambda1} shows the case 
when using potential \eqref{eq:exp_potential} with magnitude $\beta = 1$ on domain $\Omega_\epsilon = [-20, 20]$. 
The eigenvalue errors of softFEM and softIGA are also minimized when using the same softness parameters. 
The errors are $3.45\times 10^{-7}$ for softFEM and $1.66\times 10^{-9}$ for softIGA, respectively. 
This shows to a certain extent the robustness of the proposed method.

\begin{figure}
    \centering
    \includegraphics[width=\textwidth]{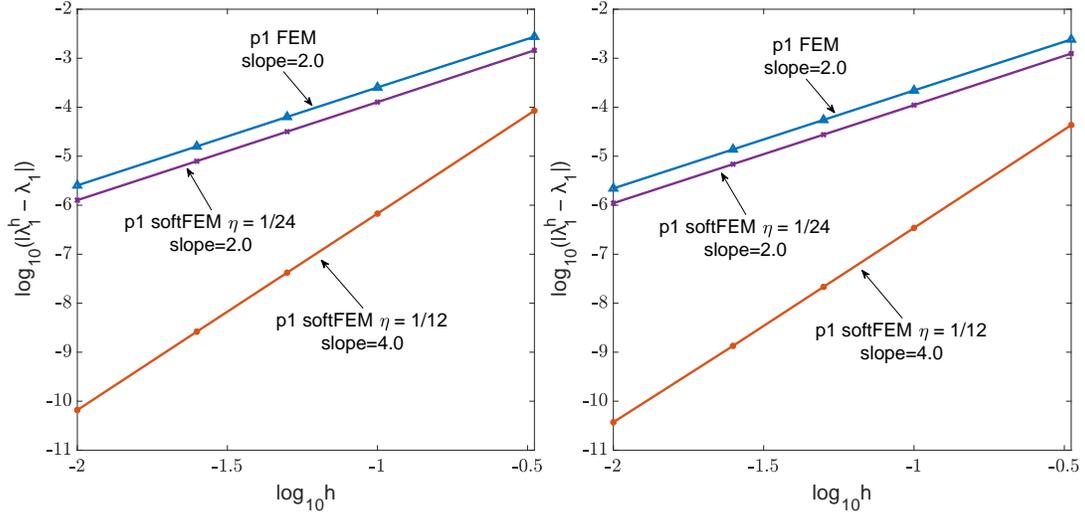}
    \caption{Eigenvalue error convergence rates of $C^0$ linear FEM and softFEM for the two-body problem with $\beta = 1$ on domain $\Omega_\epsilon = [-20,20]$. 
    The left plot shows the case of the polynomial decaying potential \eqref{eq:poly_potential} while the right plot shows that of the exponential decaying potential \eqref{eq:exp_potential}. }
    \label{fig:c0p1lambda1}
\end{figure}

\begin{figure}
    \centering
    \includegraphics[width=\textwidth]{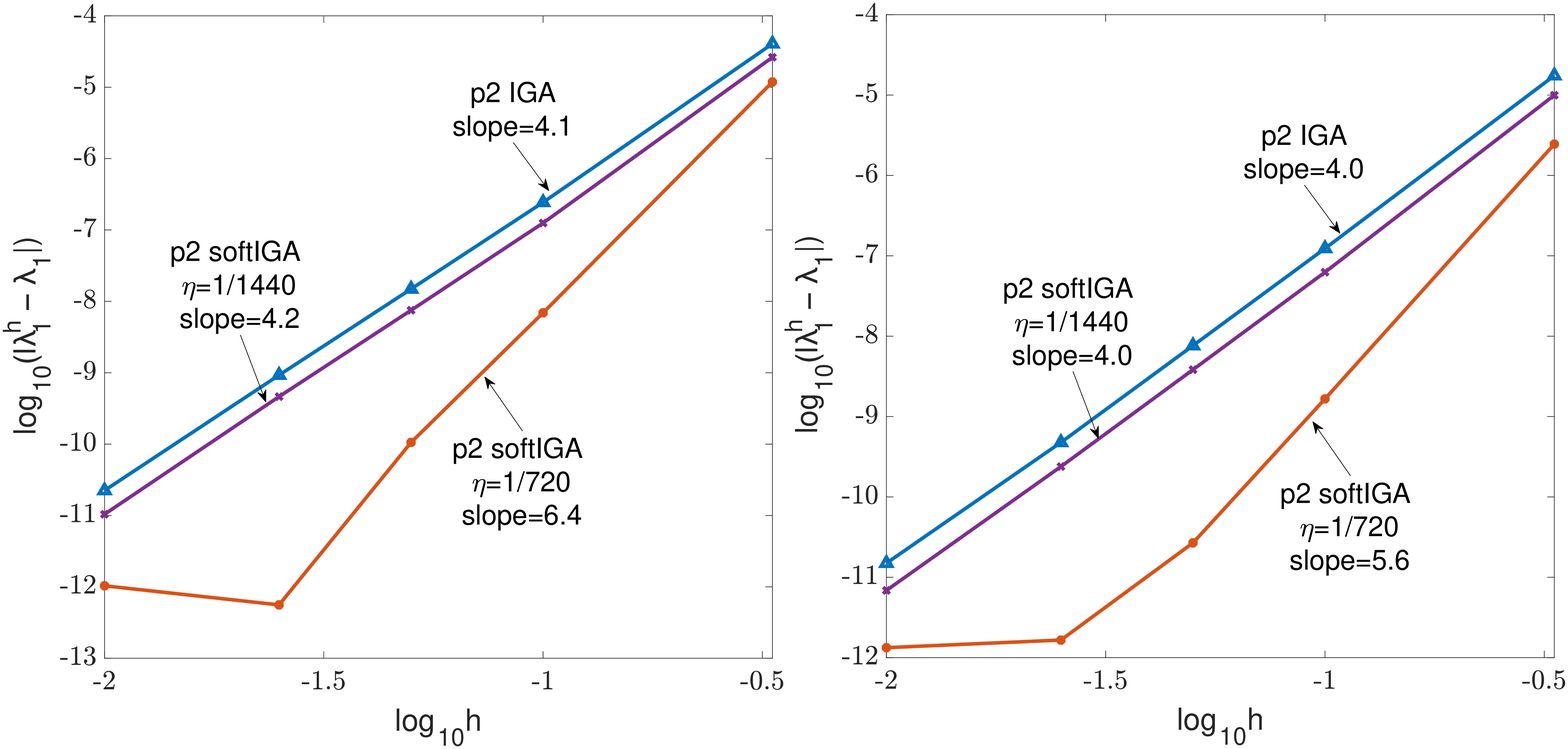}
    \caption{Eigenvalue error convergence rates of $C^1$ quadratic IGA and softIGA for the two-body problem with $\beta = 1$ on domain $\Omega_\epsilon = [-20,20]$.
    The left plot shows the case of the polynomial decaying potential \eqref{eq:poly_potential} while the right plot shows that of the exponential decaying potential \eqref{eq:exp_potential}.}
    \label{fig:c1p2lambda1}
\end{figure}

To further illustrate softFEM’s performance in one bound state, 
Figure \ref{fig:c0p1lambda1} reports the eigenvalue errors of softFEM 
with meshes of 120 to 4000 uniform elements. 
We use $\eta = 0$, $\frac{1}{24}$ and $\frac{1}{12}$, where $\eta = 0$ indicates the standard FEM.
We study both potential \eqref{eq:poly_potential} and \eqref{eq:exp_potential} with $\beta = 1$ on domain $\Omega_\epsilon = [-20, 20]$.
We use the same reference solution as above for eigenvalue errors. 
On any given mesh under consideration, we observe that softFEM with different values of $\eta$ leads to smaller eigenvalue errors than those of FEM.
The softFEM with $\eta = \frac{1}{12}$ returns the smallest eigenvalue errors around $10^{-10}$ while the rest return eigenvalue errors around $10^{-6}$ when using 4000 elements. 
The eigenvalue errors are convergent with an optimal order of $2p=2$ for linear \rev{($p=1$)} FEM and softFEM with $\eta = \frac{1}{24}$.
When $\eta=\frac{1}{12}$,
we obtain superconvergent errors of order $2p+2=4$.
Similarly, Figure \ref{fig:c1p2lambda1} shows the results when utilizing $C^1$ quadratic softIGA
for both potential \eqref{eq:poly_potential} and \eqref{eq:exp_potential}. 
With the same mesh configuration,
softIGA with $\eta = \frac{1}{720}$ returns the smallest eigenvalue errors around $10^{-12}$ while others return eigenvalue errors around $10^{-11}$.
The scale $10^{-12}$ of the errors is limited due to the accuracy of the integration involving the potential. 
SoftIGA has an optimal convergence when $\eta=0,\frac{1}{1440}$ as well as a superconvergence when $\eta = \frac{1}{720}$.

\begin{figure}
    \centering
    \includegraphics[width=\textwidth]{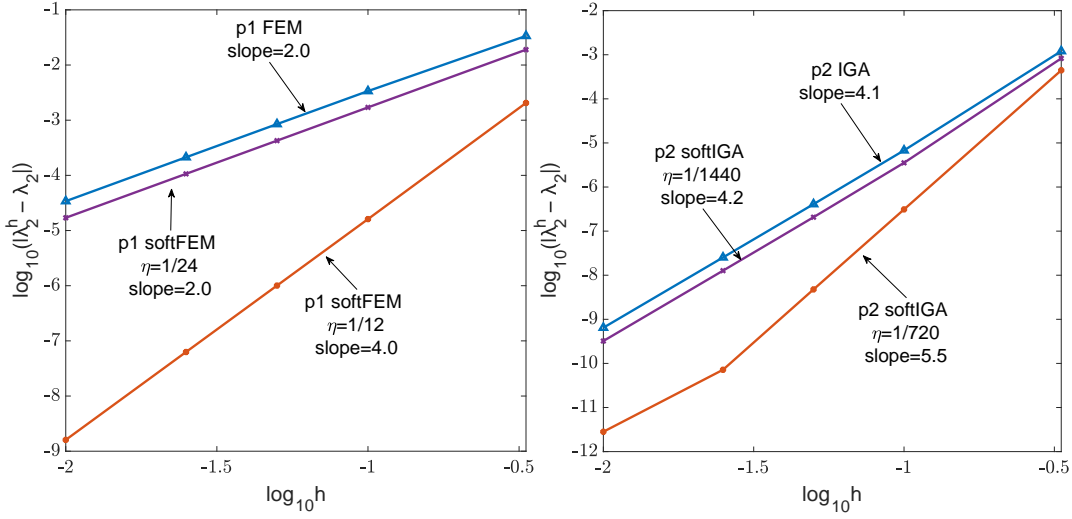}
    \caption{Eigenvalue error convergence rates of $\lambda^h_2$ for the two-body problem with $\beta = 5$ and polynomial decay \eqref{eq:poly_potential} on domain $\Omega_\epsilon = [-20,20]$.}
    \label{fig:c1p2lambda2}
\end{figure}

For a case of \eqref{eq:pde} with two bound states,
Figure \ref{fig:c1p2lambda2} shows eigenvalue errors of $\lambda_2$ 
of polynomial decay \eqref{eq:poly_potential} with $\beta = 5$.  
Again, we consider $\eta=\{0, \frac{1}{24}, \frac{1}{12}\}$ 
and $\eta=\{0, \frac{1}{1440}, \frac{1}{720}\}$ 
for softFEM and softIGA, respectively. 
With a mesh of 4000 elements,
the eigenvalue error of softFEM with $\eta=\frac{1}{12}$ reaches $10^{-9}$, 
while the errors of the other two cases are around $10^{-5}$. 
Similar behaviour is observed in the case of $C^1$ quadratic softIGA.
In summary, both softFEM and softIGA outperform their standard versions. 
This confirms our theoretical expectations.

\subsection{A Study on Domain Size}

\begin{figure}
    \centering
    \includegraphics[width=\textwidth]{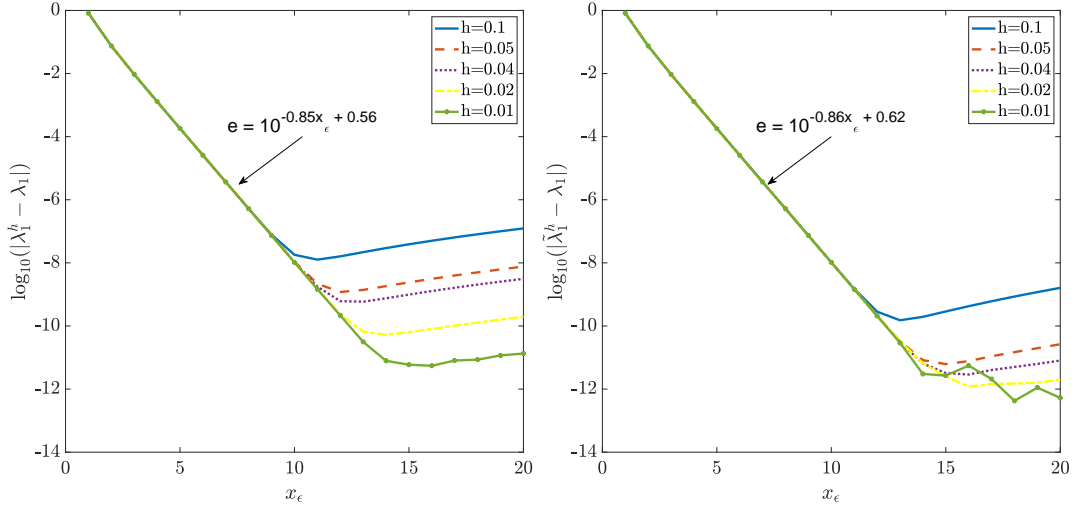}
    \caption{Eigenvalue errors with respect to domain size $x_\epsilon$ for the two-body problem with polynomial decaying potential \eqref{eq:exp_potential}. The left plot shows the case of IGA while the right plot shows that of softIGA with $\eta = 1/720$.}
    \label{fig:f1domainsize}
\end{figure}

In the previous section, we discussed the variation of eigenvalue error with $h$ within a certain interval $\Omega_\epsilon=[-20,20]$. 
The \rev{accuracy} of the problem is not only determined by the size of $h$, 
but also related to the domain size. 
To study the impact of the domain size $x_\epsilon$ on the accuracy, 
we apply $C^1$ quadratic IGA and softIGA on a one-bound state problem. 
We use the potential \eqref{eq:exp_potential} with $\kappa = 1/2$ and the potential magnitude $\beta=1$.
We apply uniform mesh grids with different sizes $h$.
Figure \ref{fig:f1domainsize} shows the relation between $x_\epsilon$ and eigenvalue errors. 
The left sub-figure is for IGA's results and the right is for SoftIGA with $\eta=1/720$.
The two fitted functions shown in Figure \ref{fig:f1domainsize} are 
$e = 10^{-0.85x_\epsilon + 0.56}$ and $e = 10^{-0.86x_\epsilon + 0.62}$, respectively. 
The smallest eigenvalue errors are $10^{-11}$ and $10^{-12}$ for IGA and softIGA, which are due to the numerical errors of the integration involving the potential function. 
The overall error is from the discretization and approximation of the domain. 
From the left sub-figure, 
we could observe the IGA discretization error dominates the overall error when the domain is large. 
Compared to the left one, 
the right sub-figure shows that softIGA reduces the IGA discretization error to some extent.
For example, if we seek an accuracy of $10^{-10}$, 
softIGA only needs a mesh with size $h=0.1$ while IGA needs a mesh with a size of 0.02 or smaller (much finer mesh) on the domain $[-13,13]$.

\subsection{Three-Body Problem with a Mass Ratio}


For solving three-body problems, we use the eigenvalues from classical BO approximation as reference solutions \cite{scherrer2017mass, happ2019universality}. 
\rev{We adopt tensor-product meshes and first} rewrite equation \eqref{eq:gmevp} as 

\begin{equation}
    (\tilde{K}^{2D} - \eta S^{2D})\tilde{U} =\tilde{ \lambda}^h M^{2D}\tilde{U}
\end{equation}
where 
\begin{equation} \label{eq:mat2d}
    \begin{array}{ll} 
    S^{2D} := S^{1D} \otimes M^{1D} + M^{1D} \otimes S^{1D}, \\
    \tilde{K}^{2D} := \tilde{K}^{1D} \otimes M^{1D} + M^{1D} \otimes \tilde{K}^{1D} + Q^{2D},\\
    M^{2D} := M^{1D} \otimes M^{1D}
    \end{array}
\end{equation}
with $Q^{2D}$ being the matrix from the inner product involving potential. \rev{Herein, $\otimes$ denotes the usual Kronecker product.}

\begin{figure}
    \centering
    \includegraphics[width=\textwidth]{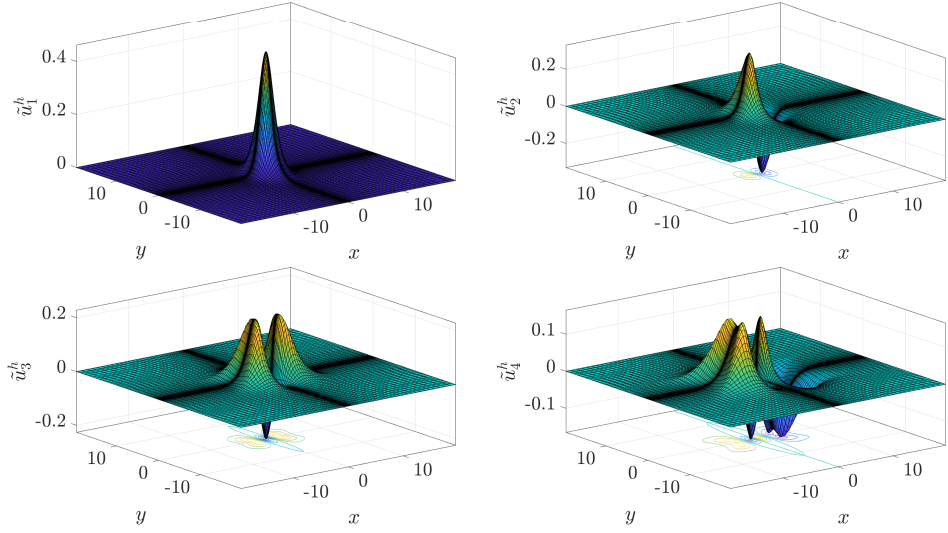}
    \caption{The first four eigenstates of the three-body problem with a mass ratio 20, potential magnitude $\beta  = 0.344595351$, and the exponential decaying potential \eqref{eq:exp_potential} using $C^0$ linear softFEM. }
    \label{fig:3bc0}
\end{figure}

\begin{figure}
    \centering
    \includegraphics[width=\textwidth]{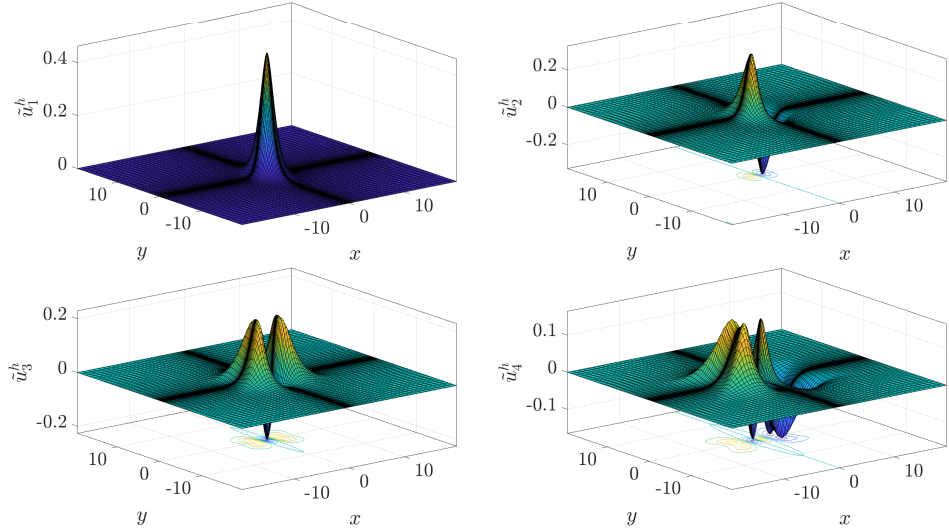}
    \caption{The first four eigenstates of the three-body problem with a mass ratio 20, potential magnitude $\beta  = 0.344595351$, and the exponential decaying potential \eqref{eq:exp_potential} using $C^1$ quadratic softIGA.}
    \label{fig:3bc1}
\end{figure}

\begin{figure}
    \centering
    \includegraphics[width=\textwidth]{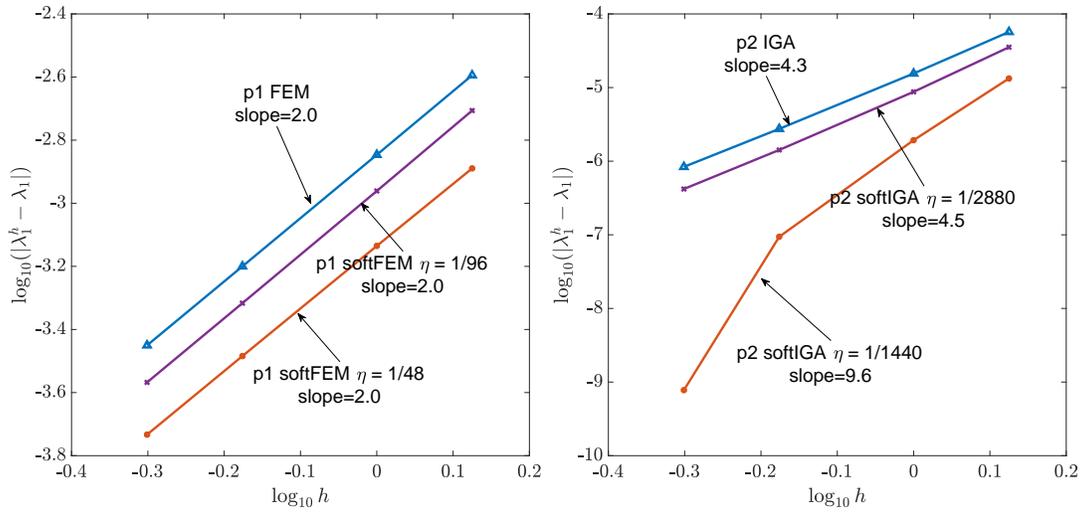}
    \caption{Eigenvalue error convergence rates of $C^0$ linear softFEM and $C^1$ quadratic softIGA for the three-body problem with potential \eqref{eq:exp_potential} and $\beta = 0.344595351$ on domain $\Omega_\epsilon = [-20,20]$.}
    \label{fig:3blambda1}
\end{figure}

For the three-body problem, we first start with a case of heavy-light mass ratio $m_h/m_l = 20$. 
$\kappa$ in the unified problem \eqref{eq:schro} is $(41/84, 0; 0, 1/21)$.
As benchmark, we use the exponential decaying potential  \eqref{eq:exp_potential} and $\beta = 0.344595351$ as the potential magnitude.
We apply $C^0$ linear softFEM elements and $C^1$ quadratic softIGA elements. 
We set the domain as $\Omega_\epsilon=[-20,20]$ and apply a non-uniform grid with $80 \times 80$ elements.
The mesh is set to be non-uniform (adaptive) for high accuracy.
Specifically, \rev{we gradually and uniformly increase} the element sizes from the origin $(0,0)$ to the boundaries such that there are $80 \times 80$ elements.
Figure \ref{fig:3bc0} shows the first four bound state eigenfunctions when using softFEM with $\eta = \frac{1}{48}$ while Figure \ref{fig:3bc1} shows the eigenfunctions of softIGA with $\eta = \frac{1}{1440}$.
These softness parameters are chosen for high accuracy as the underlying mesh is non-uniform. 
The eigenstate solution shapes match well with the ones obtained using the BO approximation in Figure 4 of \cite{happ2019universality} and Figure 7 in \cite{deng2022isogeometric}.
Similar to two-body problems, to further study the performance of softFEM and softIGA, Figure \ref{fig:3blambda1} presents the eigenvalue errors with meshes of 30 to 80 non-uniform elements. 
The eigenvalue errors of softFEM reach $10^{-3.8}$ while the errors of softIGA are around $10^{-9}$.
This figure does not show superconvergence as in two-body problems but the \rev{corresponding errors are much smaller than those of the IGA eigenvalues.
This is due to that (1) the mesh is \rev{non-uniform (the way of refining the meshes affects the convergence rate)} and (2) the potential in 2D is not separable (see \eqref{eq:mat2d}; the matrix $Q^{2D}$ is not separable; in such a case, a spatial dependent $\eta$ may lead to smaller errors).}

\subsection{Scenarios with Other Mass Ratios}

To show robustness of our method, we study the three-body problems with mass ratios $m_h/m_l = 1$ and $m_h/m_l = 100$. 
The classic BO approximation requires a large difference between electrons and nuclei's mass. 
For dynamic three-body motions with similar mass, the BO approximation is no longer applicable. 

\begin{figure}
    \centering
    \includegraphics[width=\textwidth]{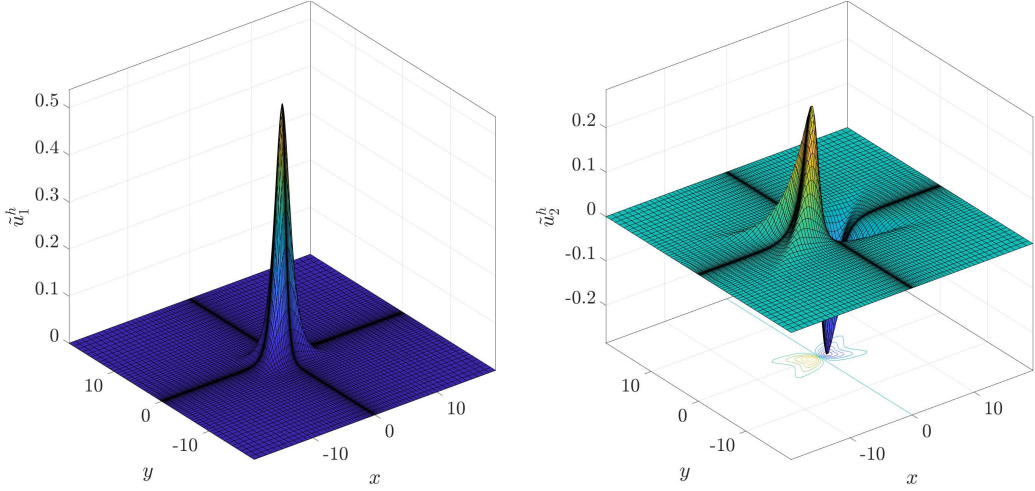}
    \caption{The first two eigenstates of the three-body problem with a mass ratio 1, potential magnitude $\beta  = 1$, and an exponential decaying potential \eqref{eq:exp_potential} using $C^1$ quadratic softIGA.}
    \label{fig:3bc1_s_1}
\end{figure}

\begin{figure}
    \centering
    \includegraphics[width=\textwidth]{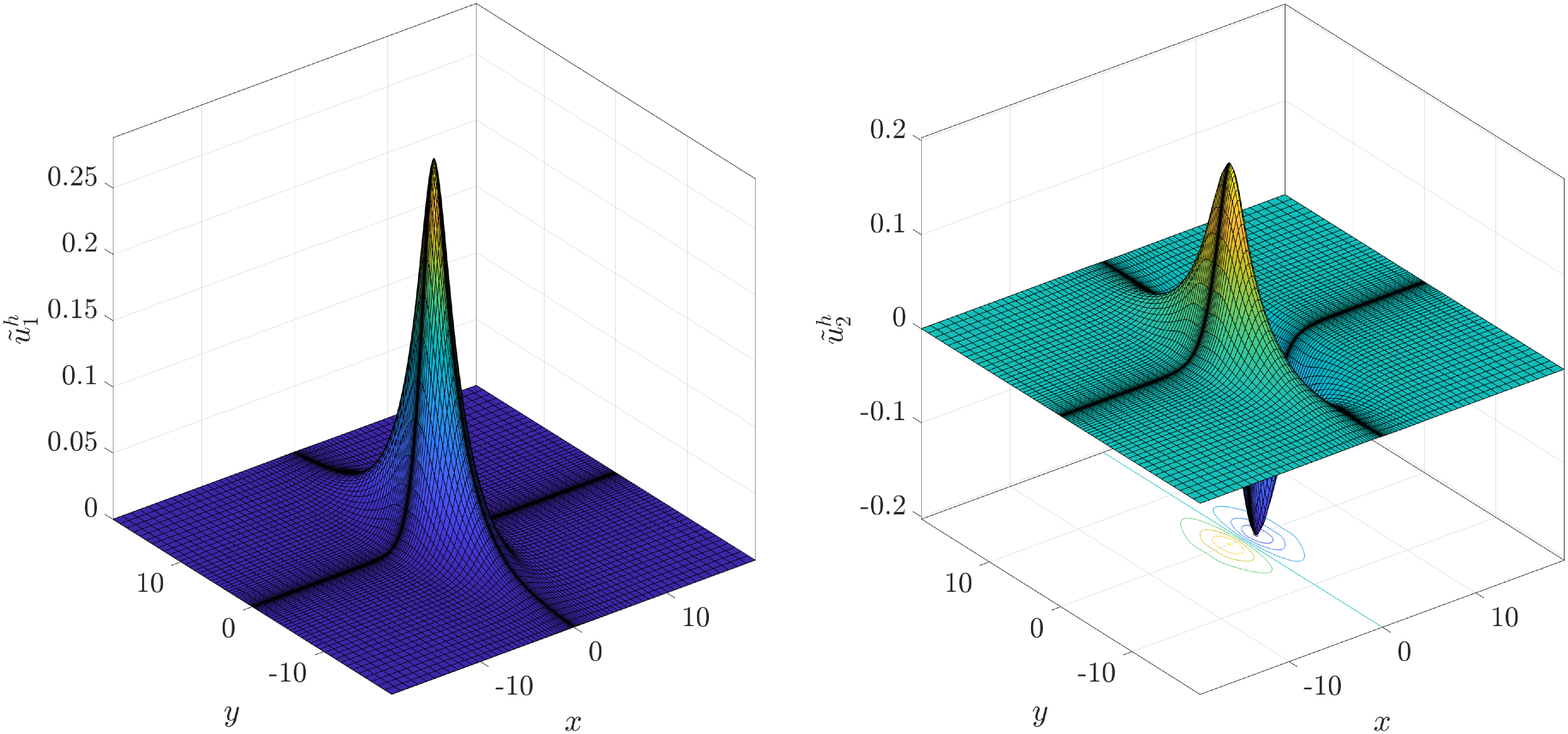}
    \caption{The first two eigenstates of the three-body problem with a mass ratio 100, potential magnitude $\beta  = 0.1$, and an exponential decaying potential \eqref{eq:exp_potential} using $C^1$ quadratic softIGA.}
    \label{fig:3bc1_s_100}
\end{figure}

First, for the unified problem \eqref{eq:schro} with a mass ratio $m_h/m_l = 1$, $\kappa=(3/8, 0; 0, 1/2)$, a potential  
\eqref{eq:exp_potential}, and a potential magnitude $\beta = 1$,
Figure \ref{fig:3bc1_s_1} shows the two bound states of $C^1$ quadratic softIGA with $\eta = \frac{1}{720}$ on $\Omega_\epsilon=[-20,20]$. 
With $80\times80$ non-uniform mesh, the eigenvalues are shown in Table \ref{tab:3beig_1}, where $\eta$ for softFEM is $\frac{1}{12}$ and $\eta$ of softIGA is $\frac{1}{720}$.
The eigenvalue errors are calculated based on reference eigenvalues obtained by using $C^6$ septic IGA with $80\times80$ non-uniform elements. 
For a case with mass ratio $m_h/m_l = 100$, $\kappa = (201/404, 0; 0, 1/101)$, 
a potential \eqref{eq:exp_potential}, and a potential magnitude $\beta = 0.1$,
Figure \ref{fig:3bc1_s_100} shows the two bound states of $C^1$ softIGA with $\eta = \frac{1}{11520}$ at $\Omega_\epsilon=[-20,20]$. 
Similarly, the eigenvalues are shown in Table \ref{tab:3beig_100}, where $\eta$ for softFEM is $\frac{1}{384}$ and $\eta$ of softIGA is $\frac{1}{11520}$.
Lastly, 
we remark that large mass ratios lead to highly heterogeneous diffusion coefficient $\kappa$,
which requires adjusting softness parameters accordingly. 
In conclusion, \rev{softFEM and softIGA elements are generally of better performance than FEM and IGA for solving the bound states of three-body problems though the outperformance can be mild in cases where mass ratio is large. }

\begin{table}[tp]
    \centering
    \begin{tabular}{|c|c|c|c|c|}
        \hline
        Method & $\lambda_1, \lambda_1^h$, or $\tilde{\lambda}_1^h$ & $\lambda_2, \lambda_2^h$, or $\tilde{\lambda}_2^h$ & Error & Error\\
		\hline 
        Reference & -0.9777963446 & -0.5425519761 & - & - \\
        \hline
        $C^0$ linear FEM & -0.9762982419 & -0.5407396655 & 1.4981e-03 & 1.8123e-03 \\
        \hline
        $C^0$ linear softFEM & -0.9780472885 & -0.5425343131 & 2.5094e-04 & 1.7663e-05 \\
        \hline
        $C^1$ quadratic IGA & -0.9777914437 & -0.5425448803 & 4.9009e-06 & 7.0958e-06 \\
        \hline
        $C^1$ quadratic softIGA & -0.9777970864 & -0.5425520381 & 7.4175e-07 & 6.1979e-08 \\
        \hline
    \end{tabular}
    \caption{Eigenvalues of the three-body problem approximated by $C^1$ quadratic IGA and softIGA with mass ratio 1 and potential \eqref{eq:exp_potential}. }
    \label{tab:3beig_1}
\end{table}

\begin{table}[tp]
    \centering
    \begin{tabular}{|c|c|c|c|c|}
        \hline
        Method & $\lambda_1, \lambda_1^h$, or $\tilde{\lambda}_1^h$ & $\lambda_2, \lambda_2^h$, or $\tilde{\lambda}_2^h$ & Error & Error\\
		\hline 
        Reference & -0.0365878475 & -0.0286770203 & - & - \\
        \hline
        $C^0$ linear FEM & -0.0365519025 & -0.0286290619 & 3.5945e-05 & 4.7958e-05 \\
        \hline
        $C^0$ linear softFEM & -0.0365608068 & -0.0286631714 & 2.7041e-05 & 1.3849e-05 \\
        \hline
        $C^1$ quadratic IGA & -0.0365877914 & -0.0286769351 & 5.6085e-08 & 8.5210e-08 \\
        \hline
        $C^1$ quadratic softIGA & -0.0365878244 & -0.0286770457 & 2.3106e-08 & 2.5395e-08 \\
        \hline
    \end{tabular}
    \caption{Eigenvalues of the three-body problem approximated by $C^1$ quadratic IGA and softIGA with mass ratio 100 and potential \eqref{eq:exp_potential}. }
    \label{tab:3beig_100}
\end{table}

\section{\rev{A Study on Computational Cost}} \label{sec:cost}

\rev{Computational resource cost mainly have two measurements, one is memory storage cost and the other is computational running time. As we are using the sparse matrix for computation, the memory cost is small. Thus, we focus on the computational time. We use Matlab default iterative solvers. 
Table \ref{tab:2b_sparse_runtime} and \ref{tab:3b_runtime} present the computational times of different methods using in two- and three-body problems. As shown in Table \ref{tab:2b_sparse_runtime}, for larger number of elements, the running time of soft method is larger than its corresponding traditional method. Yet, for the small number of elements, $n=80$ for example, the soft method is slight faster than its corresponding traditional method. This phenomenon also shown in three-body problem. In three-body problem, the size of matrix is of growth $\mathcal{O}(n^2)$: for $n$ number of elements, the dimensions of $K$, $S$ and $M$ are $\mathcal{O}(n^2)$. 
Here in Table \ref{tab:3b_runtime} we only focus on small number of elements, and in such a case, the soft method for three-body problem is faster than its corresponding traditional FEM/IGA method.
This is expected for iterative solvers as the soft methods lead to matrix problems that are of smaller condition numbers.}

\begin{table}[tp]
    \centering
    \begin{tabular}{|c|c|c|c|c|c|}
        \hline
        Method & $n=80$ & $n=400$& $n=800$ & $n=1600$ & $n=4000$\\
        \hline
        $C^0$ linear FEM & 0.031537 & 0.175362 & 0.276650 & 0.555350 & 2.319363 \\
        \hline
        $C^0$ linear softFEM & 0.030498 & 0.178336 & 0.313189 & 0.817097 & 4.061546 \\
        \hline
        $C^1$ quadratic IGA & 0.036177 & 0.155595 & 0.310764 & 0.740833 & 3.348537 \\
        \hline
        $C^1$ quadratic softIGA & 0.035118 & 0.203704 & 0.384575 & 0.918226 & 4.635101 \\
        \hline
    \end{tabular}
    \caption{Running time of two-body problem with different mass ratio $n$ $\beta = 5$ and polynomial decay \eqref{eq:poly_potential}. }
    \label{tab:2b_sparse_runtime}
\end{table}

\begin{table}[tp]
    \centering
    \begin{tabular}{|c|c|c|c|c|}
        \hline
        Method & $n=20$& $n=40$ & $n=60$ & $n=80$\\
        \hline
        $C^0$ linear FEM & 0.237728 & 4.721456 & 47.637794 & 313.714472  \\
        \hline
        $C^0$ linear softFEM & 0.237428 & 4.566356 & 43.267819 & 311.864974 \\
        \hline
        $C^1$ quadratic IGA & 0.478067 & 6.461850 & 50.374009 & 343.483040  \\
        \hline
        $C^1$ quadratic softIGA & 0.327223 & 6.354705 & 49.365816 & 335.872322 \\
        \hline
    \end{tabular}
    \caption{Running time of three-body problem with different mass ratio $n$, $\beta = 0.344595351$ and exponential decay \eqref{eq:exp_potential}. }
    \label{tab:3b_runtime}
\end{table}

\section{Concluding Remarks} \label{sec:conc}

In this paper, we study the comparison of FEM, IGA, softFEM, and softIGA for solving the quantum two- and three-body problems. 
For simplicity, we consider both $C^0$ linear and $C^1$ quadratic elements.
\rev{For small mass ratios, 
both softFEM and softIGA outperform FEM and IGA, respectively,  while for large mass ratios, the ``soft" methods have similar performance.}
Also, for two-body problems, 
we observe eigenvalue error superconvergent rates for particular choices of softness parameters. 
Lastly, we study the performance of the proposed methods for problems with different mass ratios (of scales $10^k, k =0,1,2$). 
We demonstrate that the method is robust with respect to mass ratios,
which outperforms the classical BO approximation methods \rev{for being more robustness to the mass ratios}. 

As for future work, an interesting direction is a generalization to the two- and three-body problems in multiple dimensions, 
or more generally, the $n$-body problem.
This would pose significant challenges in numerical approximation as the overall degrees of freedom in the system would increase dramatically.
One may first develop and apply techniques to reduce the model problem dimensions.
Then the proposed method may be applied as an alternative  for solving dimension-reduced problems. 


\bibliography{ref}

\begin{thebibliography}{10}

\bibitem{berezin1991schrodinger}
F.~A. Berezin and M.~A. Shubin.
\newblock {\em The {S}chr\"{o}dinger equation}, volume~66 of {\em Mathematics
  and its Applications (Soviet Series)}.
\newblock Kluwer Academic Publishers Group, Dordrecht, 1991.

\bibitem{born1985quantentheorie}
M.~Born and W.~Heisenberg.
\newblock Zur quantentheorie der molekeln.
\newblock In {\em Original Scientific Papers Wissenschaftliche
  Originalarbeiten}, pages 216--246. Springer, 1985.

\bibitem{buffa2011isogeometric}
A.~Buffa, C.~de~Falco, and G.~Sangalli.
\newblock Iso{G}eometric {A}nalysis: stable elements for the 2{D} {S}tokes
  equation.
\newblock {\em Internat. J. Numer. Methods Fluids}, 65(11-12):1407--1422, 2011.

\bibitem{cederbaum2008born}
L.~S. Cederbaum.
\newblock Born--{O}ppenheimer approximation and beyond for time-dependent
  electronic processes.
\newblock {\em The Journal of Chemical Physics}, 128(12):124101, 2008.

\bibitem{cottrell2009isogeometric}
J.~A. Cottrell, T.~J. Hughes, and Y.~Bazilevs.
\newblock {\em Isogeometric analysis: toward integration of {CAD} and {FEA}}.
\newblock John Wiley \& Sons, 2009.

\bibitem{de1978practical}
C.~de~Boor.
\newblock {\em A practical guide to splines}, volume~27 of {\em Applied
  Mathematical Sciences}.
\newblock Springer-Verlag, New York, revised edition, 2001.

\bibitem{deng2022isogeometric}
Q.~Deng.
\newblock Isogeometric analysis of bound states of a quantum three-body
  problem in 1{D}.
\newblock In {\em Computational Science -- ICCS}, pages 333--346, Cham, 2022.
  Springer International Publishing.

\bibitem{deng2023softiga}
Q.~Deng, P.~Behnoudfar, and V.~M. Calo.
\newblock Softiga: Soft isogeometric analysis.
\newblock {\em Computer Methods in Applied Mechanics and Engineering},
  403:115705, 2023.

\bibitem{deng2021outlier}
Q.~Deng and V.~M. Calo.
\newblock Outlier removal for isogeometric spectral approximation with the
  optimally-blended quadratures.
\newblock In {\em Computational science---{ICCS}. {P}art {II}}, volume 12743 of
  {\em Lecture Notes in Comput. Sci.}, pages 315--328. Springer, Cham, 2021.

\bibitem{deng2021softfem}
Q.~Deng and A.~Ern.
\newblock Soft{FEM}: revisiting the spectral finite element approximation of
  second-order elliptic operators.
\newblock {\em Comput. Math. Appl.}, 101:119--133, 2021.

\bibitem{evans2013isogeometric}
J.~A. Evans and T.~J.~R. Hughes.
\newblock Isogeometric divergence-conforming {B}-splines for the
  {D}arcy-{S}tokes-{B}rinkman equations.
\newblock {\em Math. Models Methods Appl. Sci.}, 23(4):671--741, 2013.

\bibitem{eyges1959quantum}
L.~Eyges.
\newblock Quantum-mechanical three-body problem.
\newblock {\em Physical Review}, 115(6):1643, 1959.

\bibitem{happ2019universality}
L.~Happ, M.~Zimmermann, S.~I. Betelu, W.~P. Schleich, and M.~A. Efremov.
\newblock Universality in a one-dimensional three-body system.
\newblock {\em Phys. Rev. A}, 100(1):012709, 14, 2019.

\bibitem{happ2022universality}
L.~Happ, M.~Zimmermann, and M.~A. Efremov.
\newblock Universality of excited three-body bound states in one dimension.
\newblock {\em Journal of Physics B: Atomic, Molecular and Optical Physics},
  55(1):015301, 2022.

\bibitem{hiemstra2021removal}
R.~R. Hiemstra, T.~J. Hughes, A.~Reali, and D.~Schillinger.
\newblock Removal of spurious outlier frequencies and modes from isogeometric
  discretizations of second-and fourth-order problems in one, two, and three
  dimensions.
\newblock {\em Computer Methods in Applied Mechanics and Engineering},
  387:114115, 2021.

\bibitem{hughes2005isogeometric}
T.~J.~R. Hughes, J.~A. Cottrell, and Y.~Bazilevs.
\newblock Isogeometric analysis: {CAD}, finite elements, {NURBS}, exact
  geometry and mesh refinement.
\newblock {\em Comput. Methods Appl. Mech. Engrg.}, 194(39-41):4135--4195,
  2005.

\bibitem{manni2022application}
C.~Manni, E.~Sande, and H.~Speleers.
\newblock Application of optimal spline subspaces for the removal of spurious
  outliers in isogeometric discretizations.
\newblock {\em Comput. Methods Appl. Mech. Engrg.}, 389:Paper No. 114260, 38,
  2022.

\bibitem{mitroy2013theory}
J.~Mitroy, S.~Bubin, W.~Horiuchi, Y.~Suzuki, L.~Adamowicz, W.~Cencek,
  K.~Szalewicz, J.~Komasa, D.~Blume, and K.~Varga.
\newblock Theory and application of explicitly correlated {G}aussians.
\newblock {\em Reviews of Modern Physics}, 85(2):693, 2013.

\bibitem{raynal1970transformation}
J.~Raynal and J.~R\'{e}vai.
\newblock Transformation coefficients in the hyperspherical approach to the
  three-body problem.
\newblock {\em Nuovo Cimento A (10)}, 68:612--622, 1970.

\bibitem{scherrer2017mass}
A.~Scherrer, F.~Agostini, D.~Sebastiani, E.~Gross, and R.~Vuilleumier.
\newblock On the mass of atoms in molecules: Beyond the {B}orn-{O}ppenheimer
  approximation.
\newblock {\em Physical Review X}, 7(3):031035, 2017.

\bibitem{schmid2017quantum}
E.~W. Schmid and H.~Zie{\`g}elmann.
\newblock {\em The Quantum Mechanical Three-Body Problem: Vieweg Tracts in Pure
  and Applied Physics}.
\newblock Elsevier, 2017.

\bibitem{skorniakov1957three}
G.~V. Skorniakov and K.~A. Ter-Martirosian.
\newblock Three body problem for short range forces. {I}. {S}cattering of low
  energy neutrons by deuterons.
\newblock {\em Soviet Physics. JETP}, 4:648--661, 1957.

\end{thebibliography}


\end{document}